\documentclass{article}

\usepackage{amsmath, amsthm, amssymb}
\usepackage[all]{xy}
\usepackage{path}

\newcommand{\U}{\mathfrak{U}}
\newcommand{\K}{\mathbb{K}}
\newcommand{\lt}{\mathrm{lt}}
\newcommand{\lm}{\mathrm{lm}}
\newcommand{\supp}{\mathrm{supp}}
\newcommand{\Z}{\mathbb{Z}}
\newcommand{\C}{\mathbb{C}}
\newcommand{\N}{\mathbb{N}}
\newcommand{\F}{\mathbb{F}}

\newcommand{\im}{\mathrm{Im}} 
\newtheorem{definition}{Definition}
\newtheorem{proposition}{Proposition}
\newtheorem{theorem}{Theorem}
\newtheorem{corollary}{Corollary}
\newtheorem{lemma}{Lemma}
\newtheorem{conjecture}{Conjecture}

\title{Gr\"obner basis and Anick's resolution for $\U_{\F_2}(sl^+_3)$.}
\author{Ivan Yudin\thanks{The work is supported by the FCT Grant SFRH/BPD/31788/2006. The
financial support by CMUC and FCT gratefully acknowledged.}
}

\begin{document}
\maketitle
\section{Introduction}

Despite extension groups between modules over an algebra are very easy to define
and taught nowadays in every standard course in homological algebra, it is still
to be very difficult to compute them explicitly for a given pair of modules. One
of such problems is a computation of extension groups between Weyl modules over
the Schur algebra $S(n,r)$. It was shown in the joint work~\cite{apsme} of the
author with Ana Paula Santana that this problem is closely related to the
construction of a minimal projective resolution of the trivial module $\K$ over
Kostant form $\U_\K(sl^+_n)$ of the universal enveloping algebra of the  Lie algebra $sl^+_n$. 

In this paper we compute the first three steps of a minimal projective
resolution of $\K$ in the case $p=2$ and $n=3$. For this we use Anick's
resolution constructed in~\cite{anick}. Our result depends on the knowledge of
a Gr\"obner basis for $\U_\K(sl^+_n)$. In the last section we give several
conjectures about the Gr\"obner basis for $\U_\K(sl^+_n)$. It should be noted
that with this conjectures proved it would be easy to extend the result of
Theorem~\ref{main} to the cases $p\ge 3$, $n=3$ and $p=2$, $n\ge 4$.  

In the Section~\ref{basis} we recall the definition of Gr\"obner basis and in
the Section~\ref{anickres} the construction of the Anick's resolution. 
Then we proceed with the definition of $\U_\K(sl^+_n)$ in Section~\ref{kostant}. 
The Sections~\ref{bigbasis},~\ref{smallbasis},~\ref{firststeps} contain new
results. Note that all the results of Section~\ref{bigbasis} are proved for an
arbitrary $p$ and $n$, and they will be used in the subsequent papers. 

\section{Gr\"obner basis}
\label{basis}
Let $X$ be a set. We denote by $X^*$ the set of all words with letters
in $X$. Then $X^*$ is a \emph{free monoid} generated by $X$ with the
multiplication given by concatenation of words and the unity $e$ given
by the empty word. 
There is a partial order $\prec$ on $X^*$ given by the incusion of words. Note
that $\prec$ is the coarsest partial order on $X^*$ such that $X^*$ is an
ordered monoid with $e$ the least element of $X^*$. A \emph{monoidal order} on
$X^*$ is a total order that refines $\prec$. 

Let $\K$ be a field. We denote by $\K\left\langle X^* \right\rangle$ a vector
space spanned by $X^*$. A vector space $\K\left\langle X^* \right\rangle$ is a
\emph{free associative algebra} generated by $X$. 
We will call the elements of $X^*$ \emph{monomials}, and the elements of
$\K\left\langle X^* \right\rangle$ \emph{polynomials}.
Define the support of $p\in \K\left\langle X^* \right\rangle$ to be the set of
element in $X^*$ with non-zero coefficients in $p$. If $\le$ is a monoidal order
on $X^*$ then we define the \emph{ leading monomial} $\lm(p)$ of $p\in
\K\left\langle X^* \right\rangle$ to be the maximal element of support of
$p$ with respect $\le$. Define the  \emph{leading term} $\lm(p)$ of $p$ to be the leading
monomial of $p$ with coefficient it enters in $p$. A monoidal order $\le$ on
$X^*$ can be extended to a partial order $\le$ on $\K\left\langle X^*
\right\rangle$ by the rule
\begin{align*}
	p\le q \Longleftrightarrow & \lm(p)<\lm(q)\\
	& \lt(p)=\lt(q) \mbox{ and } p-\lt(p) \le q-\lt(q).
\end{align*}
Note that in the case $\lm(p)=\lm(q)$ but $\lt(p)\not=\lt(q)$ the polynomials 
$p$ and $q$ are incompatible.

The pair $(m,f)$, where $m$ is a monomial and $f$ an element of $\K\left\langle
X^* \right\rangle$, is called a \emph{rewriting rule} if $m> f$.
Note that every element $p\in \K\left\langle X^* \right\rangle$ gives a
rewriting rule $r(p)=(\lm(p),f)$ where $f=(p-\lt(p))/\lambda$ and $\lambda$ is the
leading coefficient of $p$. 
We will say that 
$h$ is a result of application of $(m,f)$ to $g$ if there is $m'\in\supp(g)$
such that  $m'=umv$ for some $u$, $v\in X^*$, and
$h=g-\lambda m'+\lambda ufv$, where $\lambda$ is the coefficient of $m$ in
$g$. We will write in this situation $g \to_{r} h$. If $r=r(p)$ for some
$p\in \K\left\langle X^* \right\rangle$ then we write $g\to_f h$ instead of
$g \to_{r(p)} h$. Let $S$ be a collection of rewriting rules or polynomials. Then 
$g\to_S h$ denotes that there is $r\in S$ such that $g\to_r h$. Formally, $\to_S$ is a
set relation on $\K\left\langle X^* \right\rangle$. We denote by $\to_S^*$ the
reflexive and transitive closure of $\to_S$.  
An element $g$ of $\K\left\langle X^* \right\rangle$ is called
\emph{non-reducible} with respect to the set of rewriting rules or polynomials
$S$ if $g$ is a minimal element
of $\K\left\langle X^* \right\rangle$ with respect to $\to_S^*$.

\begin{definition}
Let $A$ be an algebra over a field $\K$ and $X=\left\{\, a_i \,\middle|\, i\in I
\right\}$ a set of generators of $A$. Denote by $\pi$ the canonical projection
from $\K\left\langle X^* \right\rangle$ to $A$. We say that a subset $S$ of
$\ker\left( \pi \right)$ is a \emph{Gr\"obner basis} of $\ker(\pi)$ if
$\pi$ restricted on the vector space of non-reducible elements with respect
$\left\{\, r(p) \,\middle|\, p\in S \right\}$ is an isomorphism of $\K$-vector
spaces. A Gr\"obner basis $S$ is called \emph{reduced} if elements $p\in S$
are non-reducible with respect to $S\setminus \left\{ p
\right\}$.
\end{definition}
Suppose that $\le$ is an artinian monoidal order on $X^*$, that is every
descending chain in  $X^*$ stabilizes. Let $f\in \K\left\langle X^* \right\rangle$. If
$f$ is reducible with respect to a Gr\"obner basis then there is $f_1$ such that
$f\to_S f_1$. By definition of Gr\"obner basis $f_1<f$ with respect to
the induced ordering on $\K\left\langle X^* \right\rangle$. If $f_1$ is
reducible we can find $f_2$ such that $f_1\to_S f_1$, $f_1>f_2$ and so on. Thus
we get a descending sequence $f>f_1>f_2>\dots$. As we assumed that the ordering
$\le$
is artinian this sequence have to break. Thus there is $f'$ that is
non-reducible with respect to $S$ and $f\to_S f'$. We call $f'$ the normal
form of $f$ with respect to $S$ and denote it by $NF(f,S)$. Note that the use of
the article ``the'' is justified by the fact that $f'$ is unique. In fact
suppose there are $f'$ and $f''$ such that $f\to_S f'$ and $f\to_S f''$. Then 
$f'-f''=(f'-f)+(f-f'')\in\ker(\pi)$ is an element of the kernel of the natural
projection $\pi\colon \K\left\langle X^* \right\rangle\to A$. Moreover, all
monomials in $f'-f''$ are non-reducible with respect to $S$. Since the images of
non-reducible monomials with respect to $S$ give a basis of $A$ under the map
$\pi$ it immediately follows that $f'-f''=0$.

The notion of Gr\"obner basis is closely connected with the notion of
critical pairs. We say that two monomials $m_1$, $m_2\in X^*$ \emph{overlaps} if there are
$u$, $v$, $w\in X^*$ such that $m_1= uv$ and $m_2=vw$. Note that two given
monomials can have different overlappings. To make things more convenient we define
an \emph{overlapping} as a triple $(m,m_1,m_2)$, such that there are $u$,
$v\in X^*$ such that $m=m_1v$ and $m=um_2$. 

\begin{definition}
	A \emph{critical pair} is a triple $(w,r_1,r_2)$, where $w$ is a
word	
	and $r_1=(m_1,f_1)$, $r_2=(m_2,f_2)$ are rewriting rules such that 
	there are $u$, $v\in X^*$ with the property
	$$
	w=um_1=m_2v \mbox{ or } w=um_1v = m_2.
	$$
	A word
	$w$ is called the \emph{tip} of the critical pair $(w,r_1,r_2)$.
\end{definition}

Let $(w,r_1,r_2)$ be a critical pair with $r_1$, $r_2\in S$ and $u$, $v\in X^*$
such that $w=um_1=m_2v$ (or $w=um_1v=m_2$). It is called reducible if $uf_1-f_2v
\to_S^* 0$ (respectively $uf_1v-f_2\to_S^* 0$). 
The set of rewriting rules $S$ is called complete if all critical pairs
$(w,r_1,r_2)$ with $r_1$, $r_2\in S$ are reducible.

\begin{theorem}
	\label{groebner}Suppose $\le$ is artinian monoidal ordering on
	$X^*$. A subset $S$ of $\K\left\langle X^* \right\rangle$ is a Gr\"obner
	basis of a two-sided ideal $I\subset \K\left\langle X^*
	\right\rangle$ if and only if the set of rewriting rules $\left\{\,
	r(p)
	\,\middle|\, p\in S \right\}$ is complete. 
\end{theorem}

We shall need the following proposition
\begin{proposition}
 \label{subsystem} Suppose $R$ is a complete rewriting system in
 variables $X$ and $Y$ is a subset of $X$. We denote by $R(Y)$ the
 subset of $R$ that consist from all the rules $(m,p)$ such that
 $m\in Y^*$. If for all $(m,p)\in R(Y)$ we have $p\in \K\left\langle
 Y^* \right\rangle$ then $R(Y)$ is a complete rewriting system. 
\end{proposition}
\begin{proof}
 Suppose $f\in \K\left\langle Y^*\right\rangle$ and $f\to_{R} g$
 then $f\to_{(m,p)} g$ for some $(m,p)\in R$. Since $m\preceq m'$ for
 some $m'\in \supp(f)$ and $m'\in Y^*$ we get that $(m,p)\in
 R(Y)$. By assumption of the proposition we get $p\in \K\left\langle
 Y^*\right\rangle$. Therefore $g\in \K\left\langle Y^*\right\rangle$
 and $f\to_{R(Y)} g$. Now by repetition we get that
 $f\in\K\left\langle Y^*\right\rangle$ and $f\to_{R}^* g$ implies that
 $f\to_{R(Y)}^* g$. 

 Suppose that $(w,r_1,r_2)$ is an overlap of two rules from
 $R(Y)$ and $u$, $v\in Y^*$ are such that $w=m_1v=um_2$
 ($w=um_1v=m_2$). Then $p_1v-up_2\in \K\left\langle Y^*\right\rangle$
 ($up_1v-p_2\in \K\left\langle Y^*\right\rangle$) and $p_1v-up_2\to_R
 0$ ($up_1v-p_2\to_R 0$), since $R$ is complete. But then
 $p_1v-up_2\to_{R(Y)}
 0$ ($up_1v-p_2\to_{R(Y)} 0$), which shows that $R(Y)$ is complete.
\end{proof}

\section{Anick resolution}
\label{anickres}

The anick resolution was introduced in~\cite{anick}. Let $A$ be an algebra over
a field $\K$ and $\varepsilon\colon A\to \K$ a homomorphism of algebras. Let $X =
a_1$,\dots be a set of generators of $A$ and $S\subset \K\left\langle X^*
\right\rangle$ a reduced Gr\"obner basis with respect to a monomial ordering
$\le$ on
$X^*$. For this set of data Anick constructed a free resolution of $\K$ over
$A$, which is nowadays called anick resolution. We will describe only the first
four steps of Anick's construction under additional assumption that
$\varepsilon(x)=0$ for all $x\in X$.

First we define sets $T_k$, $k=-1,0,1,2$, that will serve as bases of $A$-free
modules $P_k$. Denote by $T_{-1}$ the set $\left\{ e \right\}$ with one element
$e$ and by  $T_0$ the set $X$.   
The set $T_1$ is the set of all leading monomials in $S$.
 Denote by $\widetilde{T}_2$ the set of all possible overlaps of
 elements of $T_1$. Every element of $\widetilde{T}_2$ is a triple
 $(w,r_1,r_2)$. We say that an overlap $(w,r_1,r_2)$ is minimal if there is no
 overlap $(w',r'_1,r'_2)$ such that $w'\prec w$. Note that if an overlap
 $(w,r_1,r_2)$ is minimal then the rules $r_1$ and $r_2$ are uniquely determined
 by $w$. In fact, suppose that $(w,r_1,r_2)$ $(w,r'_1,r'_2)\in
 \widetilde{T}_2$. Then $w=m_1v=m'_1v'$. But this means that $m_1\preceq m'_1$
 or $m'_1\preceq m_1$. Since $S$ is a reduced Gr\"obner basis it follows that
 $r_1=r'_1$. Similarly  $r_2=r'_2$. Denote by $T_2$ the set of monomials in
 $X^*$ such that there is a  minimal overlap $(w,r_1,r_2)$.  Denote for $k=-1$, $0$, $1$, $2$
 by $P_k$ the $A$-linear span of $T_k$. Let $M$ be the set of all
non-reducible monomials with respect to $S$. Then for $k=-1,0,1,2$ the set
$$
N_k = \left\{\, m.t \,\middle|\,  m\in M,\ t\in T_k \right\}
$$
is the basis of $P_k$ over $\K$.

The sets $N_k$ have a full ordering induced by the ordering $\le$ on $X^*$ via
the map $m.t\mapsto mt$. We define maps $\delta_n\colon P_n\to P_n$ and
$j_n\colon P_{n-1}\to P_n$ as follows. 
\begin{align*}
	\delta_0(m.x)& := NF(mx,S).e  \\ j_0(ux.e)& := u.x\\
	\delta_1(m.t)&:= NF(mt',S).x,\mbox{ where $t=t'x$}
\end{align*} 
Now let $m\in M$ and $x\in X$. Suppose there are  $u$, $v\in M$ such that
$m=uv$ and $vx\in T_1$. Then we define $j_1(m.x)= u.vx$. Otherwise we let
$j_1(m.x)=0$. Note that $j_1$ is well-defined as $m=uv=u'v'$ would imply that
$v\preceq v'$ 
or $v'\preceq v$ and therefore $vx\preceq v'x$ or $v'x\preceq v'x$. But since
$S$ is reduced Gr\"obner basis any two different  elements of $T_1$ are
incompatible with respect to $\preceq$ (in other words $T_1$ is an anti-chain in
the Anick's terminology). 

Let $w\in T_2$ be such that $w=m_1v=um_2$ with $m_1$, $m_2\in T_1$. Define
$\delta_2(m.w)= NF(mu,S).m_2$. 

Suppose $t\in T_1$ and $m\in M$. If $m=uv$ for some $u$, $v\in M$ such that
$vt\in T_2$ then we define $j_2(m.t)=u.vt$. Note that if such $u$ and $v$ exist
then they are unique as $S$ is a reduced Gr\"obner basis. If there is no
$u$ and $v$ with the above property then we let $j_2(m.t)=0$.

Now we define homomorphisms of left $A$-modules $d_n\colon P_n\to P_{n-1}$ and
homomorphisms of $\K$-vector spaces $i_n\colon \ker(d_{n-1})\to P_n$ for $n=0,1,2$ by
induction. Since $d_n$ is a homomorphism of free $A$-modules it is
enough to define $d_n$ on the basis elements $.t$, where $t\in T_n$.
On the other hand $i_n$ is a homomorphism of $\K$ vector spaces,
moreover we do not have any convenient basis for $\ker(d_{n-1})$. We
will define $i_n$ by induction on the leading term of $f\in
\ker(d_{n-1})$.
\begin{align*}
	d_0 (.t)& := \delta_0(.t)\\
	i_0 (m.e)& := j_0(m.e)\\
	d_{n+1} (.t)&:= \delta_{n+1}(.t) - i_nd_n(\delta_{n+1}(.t)) \\
	i_n(f) &:= j_n(\lt(f)) + i_n(f-d_n(j_n(\lt(f)))).
\end{align*}
 Note that it is not obvious that $d_n$ and $i_n$ are
 well-defined. This a part of the claim of Proposition~\ref{anick}.
The following proposition is proved in~\cite{anick}. Note that Anick~\cite{anick}
constructed modules $P_n$ and maps $d_n$ for all $n\in \N$. 
\begin{proposition}
	\label{anick} The sequence of left $A$-modules 
	$$
	P_2 \stackrel{d_2}{\longrightarrow} P_1 \stackrel{d_2}{\longrightarrow}
	P_0 \stackrel{d_2}{\longrightarrow} P{-1}   \stackrel{\varepsilon}{\longrightarrow}\K \to 0
	$$
	is an exact complex. 
\end{proposition}

\section{Konstant form of universal enveloping algebra}
\label{kostant}
Denote by $sl_n^+$ the Lie algebra of upper triangular nilpotent
matrices.
Let $\U_n(\C)$ be its universal enveloping algebra over $\C$. 
 We shall consider $sl^{+}_n$ with the standard basis $\left\{
e_{ij}\mid 1\le i<j\le n\right\}$. Then $\U_n(\C)$ is generated as an algebra by the
elements $e_{1,2},e_{2,3},\dots, e_{n-1,n} $. 

Let $\ll$ be an arbitrary full ordering on the set $\left\{\, e_{ij} \,\middle|\,
i<j \right\}$. 
We always
assume that in the product $\prod\limits_{i<j}e_{ij}^{k_{ij}}$ the
generators increase from the left to right, with respect to the 
ordering $\ll$. It follows from the Poincare-Birkhoff-Witt Theorem, that 
the set 
$$
\mathbb{B}_n =\left\{ \prod_{1\le i<j\le n} e_{ij}^{k_{ij}}\middle| k_{ij}\in \N
\right\}
$$
is a $\C$-basis of $\U_n(\C)$.
 Denote by $e_{ij}^{(k)}$
the element $\frac{1}{k!}e_{ij}^{k}$ of the algebra $\U_n(\C)$.
We define $\U_n(\Z)$ to be the $Z$-sublattice of $\U_n(\C)$ generated by
the set 
$$
\overline{\mathbb{B}}_n =\left\{ \prod_{i<j} e_{ij}^{(k_{ij})}\middle| k_{ij}\in \N
\right\}.
$$
\begin{proposition}
 The set $\U_n(\Z)$ is a subring of $\U_n(\C)$. In other words,
 $\U_n(\Z)$  is a $\Z$-algebra. It is called the \emph{Kostant form}
 of the universal enveloping algebra $\U_n(\C)$ over $\Z$. 
\end{proposition}
\begin{proof}
 For a proof see~\cite[Lemma~2 after Proposition~3]{kostant}
 and~\cite[Remark~3]{kostant} thereafter. 
\end{proof} 

\begin{definition}
 For any field $\K$, the algebra $\U_n(\K) := \K\otimes_{\Z} \U_n(\Z)$ is called
 \emph{Kostant form} of the algebra $\U_n(\C)$ over 
 $\K$.
\end{definition}
Define a grading on $sl_n^+\left( \C \right)$ by $\deg(e_{ij})= j-i$. 

 This grading extends to the grading of $\U_\C(sl^+_n)$
such that 
$$
\U_\C(sl^+_n)_{d} = \left\{\, \prod_{i<j} e_{ij}^{k_{ij}}\,\middle|\,
\sum_{i<j} k_{ij} (j-i)=d \right\}.
$$
Since the intersection of $\U_\C(sl^+_n)_d$ with $\U_\Z(sl^+_n)$ is a
lattice in $\U_\C(sl^+_n)_d$, this grading downgrades to $\U_\Z(sl^+_n)$.
After tensoring with $\K$ we get a grading on $\U_\K(sl_n^+)$ such that
$$
\U_\K(sl^+_n)_{d} = \left\{\, \prod_{i<j} e_{ij}^{(k_{ij})}\,\middle|\,
\sum_{i<j} k_{ij} (j-i)=d \right\}.
$$

\section{Big Gr\"obner basis}
\label{bigbasis}
In this section we describe a 
Gr\"obner basis of the algebra $\U_n(\K)$ with respect to the generating set
$X =\left\{\, e^{(k)}_{ij} \,\middle|\, i<j,\ k\in \N \right\}$. 
We will consider deglex ordering on $X^*$ with respect to the degree function defined in the previous section and the ordering $\ll$
on $X$. 

\begin{theorem}
Let $X$ and the ordering on $X$ be as above. Then the following set of rewriting rules is complete:
	\begin{align}
		e_{ij}^{(k)} e_{ij}^{(r)} &\to {k+r \choose k}
		e_{ij}^{(k+r)}\label{one}\\
		e_{ij}^{(k)} e_{st}^{(r)} &\to e_{st}^{(r)} e_{ij}^{(k)}\mbox{,
		if $i$, $j$, $s$, $t$ are different and $(s,t)\ll
		(i,j)$}\label{two}\\
		e_{ij}^{(k)} e_{jt}^{(r)} &\to \sum_{s=0}^{min(k,r)}
		e_{jt}^{(r-s)} e_{i,t}^{(s)} e_{i,j}^{(k-s)}\mbox{, if $(i,j)\gg
		(j,t)$}\label{three}\\
e_{ij}^{(k)} e_{si}^{(r)} &\to \sum_{t=0}^{min(k,r)}(-1)^t e_{si}^{(r-t)} e_{s,j}^{(t)} e_{ij}^{(k-t)}\mbox{, if $(i,j)\gg (s,i)$}
\label{four}
	\end{align}
	\label{big}
\end{theorem}
Note that the corresponding Gr\"obner basis is not reduced in general, since it
can happen that $(i,t)$ doesn't lie between $(j,t)$ and $(i,j)$ in~(\ref{three})
or that $(s,j)$ doesn't lie between $(s,i)$ and $(i,j)$ in~(\ref{four}). 
\begin{proof}
	It is clear that the set 
	$$
	B = \left\{\, \prod_{i<j}e_{i,j}^{(k_{ij})} \,\middle|\, i<j,\ k_{ij}\in \N \right\}
	$$
	is the set of non-reducible words with respect to the given rewriting
	system. By definition the natural image of $B$ in
	$\U_{\K}(sl^+_n)$ is a basis of $U_{\K}(sl_n^+)$. Therefore it is enough
	to check that for every rule the left hand side and the right hand side
	are equal in $\U_\K(sl^+_n)$. 
	This is obvious for (\ref{one}) and (\ref{two}). Thus we have only to
	check the claim for (\ref{three}) and (\ref{four}). We shall do this
	only for (\ref{three}) as the case (\ref{four}) is similar. 

	We have to prove the equality
	$$
	e_{ij}^{(k)} e_{jt}^{(r)} \to \sum_{s=0}^{min(k,r)}
		e_{jt}^{(r-s)} e_{i,t}^{(s)} e_{i,j}^{(k-s)}
		$$
	in $\U_\K(sl^+_n)$. Clearly it is enough to prove the same equality in
	$\U_\Z(sl^+_n)$ and therefore in $\U_\C(sl^+_n)$. We will do this by
	induction on the minimum of  $k$ and $r$. The case $min(k,r)=1$ splits
	into two cases $r=1$ and $k=1$. The case $k=1$ we prove by induction on
	$k$. For $k=r=1$ we have
	$$
	e_{ij} e_{jt}= e_{jt}e_{ij} + e_{it}.
	$$
	Suppose we have proved equality for $k=1$ and $r\le r_0$. Let us check
	it for $r=r_0+1$. 
	\begin{align*}
		e_{ij}e_{jt}^{(r)} &=\frac{1}{r} 
		e_{ij}e_{jt}^{(r-1)}e_{jt} \mbox{\ \ \ \ \ \ \ \ \ \     induction assumption} \\
		&= \frac{1}{r} \left(e_{jt}^{(r-1)}e_{ij}+
		e_{jt}^{(r-2)}e_{it}\right)e_{jt} \\
		&= \frac{1}{r}\left( e_{ij}^{(r-1)}e_{jt}e_{ij} +
		e_{jt}^{(r-1)} e_{it} + e_{jt}^{(r-2)}e_{jt}e_{it}\right) \\
		&= e_{ij}^{(r)} e_{ij} + \frac{1}{r}\left( 1 + r-1)
		\right)e_{jt}^{(r-1)} e_{it} \\
&= e_{ij}^{(r)} e_{ij} + e_{jt}^{(r-1)} e_{it}. 
	\end{align*}
	Now we prove the equality in the case $r=1$ and $k\ge 2$. Suppose it is
	proved for all $k\le k_0$. Let us show it for $k=k_0+1$. We have
	\begin{align*}
		e_{ij}^{(k)} e_{jt} & =
		\frac{1}{k}e_{ij}e_{ij}^{(k-1)}e_{jt} \\
		&= \frac{1}{k} \left( e_{ij} e_{jt} e_{ij}^{\left( k-1
		\right)}+ e_{ij} e_{it}e_{ij}^{(k-2)}  \right)\\
		&= \frac{1}{k} \left( e_{jt}e_{ij}e_{ij}^{(k-1)}
		+e_{it}e_{ij}^{(k-1)} + e_{it}e_{ij}e_{ij}^{(k-2)}  \right)\\
		&= e_{jt}e_{ij}^{(k)} + e_{it} e_{ij}^{(k-1)}.
	\end{align*}
Suppose we have prove equality for all $k$ and $r$ such that $min(k,r)\le m_0$.
Let us prove it for $min(k,r)= m_0+1$. There are two cases $k=m_0+1$ and
$r=m_0+1$. As the computations are very similar we will treat only the first
case. 
\begin{align*}
	e_{ij}^{(k)}e_{jt}^{(r)}&=\frac{1}{k}e_{ij}e_{ij}^{(k-1)}
	e_{jt}^{(r)}\\
	&= \frac{1}{k}\sum_{s=0}^{k-1} e_{ij}e_{jt}^{(r-s)}
	e_{it}^{(s)}e_{ij}^{k-1-s} \\
	&= \frac{1}{k} \sum_{s=0}^{k-1}\left( e_{jt}^{(r-s)} e_{it}^{(s)} e_{ij}
	e_{ij}^{(k-1-s)} + e_{(jt)}^{(r-s-1)} e_{it}e_{it}^{(s+1)}
	e_{ij}^{(k-1-s)}\right) \\
	&= \frac{1}{k} \sum_{s=0}^k (k-s + s) e_{jt}^{(r-s)} e_{it}^{(s)}
	e_{ij}^{(k)}= \sum_{s=0}^k e_{jt}^{(r-s)} e_{it}^{(s)}
	e_{ij}^{(k)}.
\end{align*}
\end{proof}

\begin{corollary}
	\label{subalgebra}
	Let $p$ be a characteristic of the field $\K$ and $l\ge 0$. Then the
	linear span $\U_\K^l(sl^+_n)$ of the set
	$$
	B'= \left\{\, \prod_{i<j} e_{ij}^{(k_{ij})} \,\middle|\, k_{ij}\le p^l-1 \right\}
	$$
	is a subalgebra of $\U_\K(sl_n^+)$
\end{corollary}
\begin{proof}
We claim that $\U_\K^l(sl^+_n)$ is the subalgebra $A$  of $\U_\K(sl^+_n)$ generated
by the set 
$$
X'= \left\{\, e_{ij}^{(k)} \,\middle|\,  k\le p^l-1 \right\}. 
$$
It is enough to show that the set $B'$ is a basis of $A$.  
Let $R$ be rewriting system defined in Theorem~\ref{big}. We claim
that $R(X')$ is complete. To prove this we apply
Proposition~\ref{subsystem}. It is obvious for the rules
(\ref{two}), (\ref{three}) and (\ref{four}) that if the left hand side
is an element of $\left( X' \right)^*$ then all the monomials on the
right hand side are also elements of $\left( X' \right)^*$. Moreover,
if $k+r\le p^l-1$ then the same is true for the rewriting rule
(\ref{one}). Suppose $k$, $r\le p^l-1$ and $k+r\ge p^l$. Then
${k+r choose k}=0$ in $\K$. It is well known that the degree of
$p$ in the prime decomposition of $n!$ is given by the formula
$$
\sum_{j=0}^{\infty}\left[ \frac{n}{p^j} \right].
$$
Therefore the degree of $p$ in the prime decomposition of ${k+r
\choose k}$ is 
\begin{align*}
 \sum_{j=0}^{\infty}\left( \left[ \frac{k+r}{p^j} \right]-\left[
 \frac{k}{p^j} 
 \right]-\left[ \frac{r}{p^j} \right] \right) &=
 \left[ \frac{k+r}{p^l} \right] + \sum_{j=0}^{l-1}
\left( \left[ \frac{k+r}{p^j} \right]-\left[
 \frac{k}{p^j} 
 \right]-\left[ \frac{r}{p^j} \right] \right)\\
 &\ge \left[ \frac{k+r}{p^l} \right]=1>0.
\end{align*}
Therefore for the rule (\ref{one}) and $k+r\ge p^l$ we get
$$
e_{ij}^{(k)}e_{ij}^{(r)} \to 0.
$$
This shows that $R(X')$ is complete. Now it is obvious that $B'$ is
the set of non-reducible monomials in the alphabet $X'$ with respect to
the rewriting system $R(X')$. This shows that $B'$ is a basis of the
algebra $A'$.
\end{proof}

\section{Small Gr\"obner basis}
\label{smallbasis}
The Gr\"obner basis obtained in the previous section is not convenient
for the construction of minimal projective resolution of $\K$, since
the anick resolution is much closer to the minimal resolution if the
chosen generating set is minimal. 

In this section we will stick to the case $p=2$ and $n=3$. The more
general case will be considered in other paper. Nevertheless we start
with the technical result that is true for an arbitrary $p$ and
$n$.
\begin{lemma}
 \label{powers} Let $p$ be the characteristic of $\K$ and
 $k=k_lp^l+k_{l-1}p^{l-1}+\dots + k_0$ with $0\le k_s\le p-1$. Then 
 for any $i<j$
 $$
 \prod_{s=0}^l \left(e_{ij}^{(p^s)}\right)^{k_s}
 $$
 is a non-zero multiple of $e_{ij}^{(k)}$.
\end{lemma}
\begin{proof}
 We have to check that the integer 
$$
n := \frac{k!}{\prod_{s=0}^l \left( p^s! \right)^{k_s}}
$$
is non-zero in $\K$. The degree of $p$ in the prime decomposition of
$n$ is given by
\begin{align*}
 \sum_{t=0}^l \left( \left[ \frac{k}{p^t} \right]-\sum_{s=0}^l
 k_s\left[ \frac{p^s}{p^t} \right] \right) & =
 \sum_{t=0}^l \left(\left( \sum_{s=t}^l  k_sp^{s-t} \right) -
 \sum_{s=t}^l k_s p^{s-t} \right) = 0.
\end{align*}
This shows that $n$ is non-zero in $\K$.
\end{proof}
Now we note that 
$$
e_{ij}^{(p^l)} = e_{i,i+1}^{(p^l)}e_{i+1,j}^{(p^l)}-
\sum_{s=0}^{p^l-1} e_{i+1,j}^{(p^l-s)} e_{ij}^{(s)}
e_{i,i+1}^{(p^l-s)}.
$$
In fact it was proved in Theorem~\ref{big}. From this equality by
induction on $j-i$ and $l$ it follows that the set 
$$
\left\{ e_{i,i+1}^{(p^l)}\,\middle|\,\, 1\le i\le n-1,\ l\in\N_0 \right\}
$$
generates $\U_\K\left( sl^+_n \right)$. 
Note that the set 
$$
\left\{ e_{i,i+1}^{(p^k)}\,\middle|\,\, 1\le i\le n-1,\ 0\le k\le l \right\}
$$
generates the subalgebra $\U_\K(sl^+_n)$ of $\U_\K(sl^+_n)$.

From now on we assume that $n=3$ and $p=2$. For a convenience we will denote
$e_{12}^{2^k}$ by $a_k$ and $e_{23}^{2^k}$ by $b_k$. 
We start with the proof of some equalities in $\U_\K(sl^+_3)$.
\begin{proposition}
	\label{squares}
	For any $k$ we have $a_k^2 = b_k^2 =0$. 
\end{proposition}
\begin{proof}
	In the proof of Corollary~\ref{subalgebra} it was proved that if
	$r$, $s\le p^t-1$ and $r+s\ge p^t$ then for any $i<j$ we have
	$e_{ij}^{(r)}e_{ij}^{(s)}=0$. We apply this claim to the situation
	$(i,j)=(1,2)$, $(2,3)$ and $r=s=2^k$, $t=k+1$.
\end{proof}
\begin{proposition}
	\label{commute}
	For any $l$ and $k$ elements $a_l$ and $a_k$ commutes. Similarly
	$b_l$ and $b_k$.
\end{proposition}
\begin{proof}
	Obvious.
\end{proof}
\begin{proposition}
	\label{skew}
	For any $l>k$ we have
	\begin{align}
		a_lb_k& + b_ka_l + a_kb_ka_ka_{k+1}\dots a_{l-1} =0\\
		b_la_k&+ a_kb_l + b_ka_kb_kb_{k+1}\dots b_{l-1}=0
		\end{align}
		in $\U_\K(sl^+_n)$.
\end{proposition}
\begin{proof}
	 We have
	\begin{align*}
		a_lb_k &= e_{12}^{(2^l)} e_{23}^{(2^k)} \\
		&= \sum_{s=0}^{2^k} e_{23}^{(2^k-s)} e_{13}^{(s)}
		e_{12}^{(2^l-s)}& \mbox{ relation (\ref{three})}\\
		&= b_ka_l + \sum_{s=1}^{2^k} e_{23}^{(2^k-s)} e_{13}^{(s)}
		e_{12}^{(2^l-s)} 
	\end{align*}
	and
	\begin{align*}
		a_kb_ka_k\dots a_{l-1} &= e_{12}^{(2^k)}
		e_{23}^{(2^k)}e_{12}^{(2^k)} \dots e_{12}^{(2^{l-1})} \\
		&= e_{12}^{(2^k)}
		e_{23}^{(2^k)}e_{12}^{(2^{l}-2^k)}& \mbox{  Lemma~\ref{powers}} \\
		&= \sum_{s=0}^{2^k} e_{23}^{(2^k-s)} e_{13}^{(s)}
		e_{12}^{(2^k-s)} e_{12}^{(2^l-2^k)}. 
	\end{align*}
	Since $2^k$, $2^l-2^k\le 2^l-1$ and $2^l=2^k+(2^l-2^k)\ge 2^l$ by
	proof in Corollary~\ref{subalgebra} we have $e_{12}^{(2^k)}
	e_{12}^{(2^l-2^k)}=0$. Now if $1\le s \le 2^k$ then $2^k-s =
	2^{i_1} + \dots + 2^{i_\sigma}$ for some $0\le \sigma\le k$ and
	$0\le i_1<i_2<\dots <i_{\sigma}\le k-1$. Moreover $2^l-2^k =
	2^k+2^{k+1}+\dots+2^{l-1}$. Therefore applying Lemma~\ref{powers} twice
	we get for $1\le s\le 2^k$
	$$
	e_{12}^{(2^k-s)} e_{12}^{(2^l-2^k)}=e_{12}^{(2^l-s)}.
	$$
	Therefore
	$$
		a_kb_ka_k\dots a_{l-1} =  \sum_{s=1}^{2^k} e_{23}^{(2^k-s)} e_{13}^{(s)}
		e_{12}^{(2^l-s)}.
	$$
	The second equality follows from the obvious duality $a_k
	\leftrightarrow b_k$.
\end{proof}

\begin{proposition}
	\label{braid}
For all $k\in \N$ we have $b_ka_kb_ka_k + a_kb_ka_kb_k=0$ in  $\U_\K(sl^+_n)$.
\end{proposition}
\begin{proof}
	We have
	\begin{align*}
		a_kb_ka_kb_k &= e_{12}^{(2^k)} e_{23}^{(2^k)}e_{12}^{(2^k)}
		e_{23}^{(2^k)}\\
		&= \sum_{s=0}^{2^k} e_{12}^{(2^k)} e_{12}^{(s)}
		e_{13}^{(2^k-s)} e_{23}^{(s)}e_{23}^{(2^k)} & \mbox{relation
		(\ref{four})}.
	\end{align*}
	Now if $0\le s \le 2^k-1$ applying Lemma~\ref{powers} twice we get
	$e_{12}^{(2^k)} e_{12}^{(s)} = e_{12}^{(2^k_s)}$. Similarly for
	$e_{13}$. If $s=2^k$ we get $e_{12}^{(2^k)}e_{12}^{(s)}=a_k^2=0$ by
	Proposition~\ref{squares}.
Therefore
\begin{align*}
		a_kb_ka_kb_k &= \sum_{s=0}^{2^k-1} e_{12}^{(2^k+s)}
		e_{13}^{(2^k-s)} e_{23}^{(2^k+s)}.
\end{align*}
From the duality $e_{12}\leftrightarrow e_{23}$ it follows that 
\begin{align*}
		b_ka_kb_ka_k &= \sum_{t=0}^{2^k-1} e_{23}^{(2^k+t)}
		e_{13}^{(2^k-t)} e_{12}^{(2^k+t)}.
\end{align*}
Now we can use rewriting rules (\ref{two}), (\ref{four}), and (\ref{one}):
\begin{align*}
	\sum_{t=0}^{2^k-1} e_{23}^{(2^k+t)}
	e_{13}^{(2^k-t)} e_{12}^{(2^k+t)} &= \sum_{t=0}^{2^k-1}
	\sum_{r=0}^{2^k+t} e_{12}^{(2^k+t-r)} e_{13}^{(r)} e_{13}^{(2^k-t)}
	e_{23}^{(2^k+t-r)} \\
	&= \sum_{s=-2^k}^{2^k-1} \sum_{t=max(0,s)}^{2^k-1}{2^k-s \choose 2^k-t}
	e_{12}^{(2^k+s)} e_{13}^{(2^k-s)} e_{23}^{(2^k+s)}.
\end{align*}
Suppose $-2^k\le s\le -1$. Denote $-s$ by $\widetilde{s}$. Then
\begin{align*}
	 \sum_{t=max(0,s)}^{2^k-1}{2^k-s \choose 2^k-t}&=
	 \sum_{t=0}^{2^k-1}{2^k+\widetilde{s} \choose 2^k-t}=
	 \sum_{j=1}^{2^k}{2^k+\widetilde{s} \choose
	 j}=1 + \sum_{j=0}^{2^k}{2^k+\widetilde{s} \choose
	 j}.
\end{align*}
Now the sum $\sum_{j=0}^{2^k}{2^k+\widetilde{s} \choose
j}$ is a coefficient of $x^{2^k}$ in the product 
	 \begin{align*}
		 \left( \sum_{j=0}^{2^k+\widetilde{s}}x^j \right) \left(
		 \sum_{j=0}^{\infty}x^j
		 \right) &= (1+x)^{2^k+\widetilde{s}} (1+x)^{-1} =
		 (1+x)^{2^k+\widetilde{s}-1}\\
		 &= (1+x)^{2^k} (1+x)^{\widetilde{s}-1} = (1 +
		 x^{2^k})(1+x)^{\widetilde{s}-1}.
	 \end{align*}
	 Since $0\le \widetilde{s}-1\le 2^k-1$, this coefficient is $1$.
	 Therefore 
\begin{align*}
	 \sum_{t=max(0,s)}^{2^k-1}{2^k-s \choose 2^k-t}&=0
\end{align*}
for $-2^k\le s \le -1$. Suppose $0\le s \le 2^k-1$. Then we get
\begin{align*}
	 \sum_{t=max(0,s)}^{2^k-1}{2^k-s \choose 2^k-t}&=
\sum_{t=s}^{2^k-1}{2^k-s \choose 2^k-t}=\sum_{j=1}^{2^k-s}{2^k-s \choose
j}\\& =1+ (1+1)^{2^k-s} = 1. 
\end{align*}
Therefore 
\begin{align*}
		b_ka_kb_ka_k &=  \sum_{s=0}^{2^k-1} e_{12}^{(2^k+s)}
		e_{13}^{(2^k-s)} e_{23}^{(2^k+s)}
\end{align*}
as required.
\end{proof}
	Let $X_l=\left\{\, a_k, b_k \,\middle|\, 0\le k\le l \right\}$.
	We order the elements of $X_l$ by
	$$
	a_0 < b_0 < a_1 < b_1 <\dots < a_l < b_l
	$$
	and define degree $\deg\colon  X_l\to \N$ by $\deg(a_k)=\deg(b_k)=2^k$.
	Denote by $\pi$ the natural projection $\K\left\langle X^*_l\
	\right\rangle\to \U_\K^l\left( sl_n^+ \right))$.
\begin{proposition}
	\label{gbl}
	The following set $S_l$ of elements in $\K\left\langle X_l^* \right\rangle$
	\begin{align}
		a_lb_k& + b_ka_l + a_kb_ka_ka_{k+1}\dots a_{l-1}& \mbox{ if
		$l>k$}\\
		b_la_k&+ a_kb_l + b_ka_kb_kb_{k+1}\dots b_{l-1}&
		\mbox{ if $l>k$}\\
		a_la_k& + a_ka_l\mbox{ if $l>k$}\\
		b_lb_k & + b_kb_l\mbox{ if $l>k$}\\
		b_ka_kb_ka_k&+  a_kb_ka_kb_k\\
		a_k^2\\
		b_k^2
	\end{align}
	is a reduced Gr\"obner basis of $\ker(\pi)$. 
\end{proposition}
\begin{proof}
	If follows from Propositions~\ref{squares},~\ref{commute},~\ref{skew},~\ref{braid}, that
	$S_l$ is a subset of $\ker(\pi)$. 
	Thus it is enough to show that the images of non-reducible monomials in
	$X_l^*$ give a basis of $\U^l_\K(sl^+_n)$. Since the images of
	non-reducible monomials in 
	$X^*_l$ generate $\U^l_\K(sl^+_3)$ as a vector space and
	$\U^l_\K(sl^+_3)$ is finite dimensional, it is enough to show that
	the number of non-reducible monomials in $X_l^*$ with respect to $S_l$ is
	less or equal to the dimension of $U^l_\K(sl^+_3)$. From
	Corollary~\ref{subalgebra} it follows that the dimension of
	$\U_\K^l(sl^+_3)$ is $\left( 2^l \right)^3= 2^{3l}=8^l$. 

	To find non-reducible monomials with respect to $S_l$ it is enough to find
	monomials that does not contain submonomials $a_lb_k$, $b_la_k$, $a_la_k$,
	$b_lb_k$ for $l>k$, and submonomials $b_ka_kb_ka_k$, $a_k^2$, $b_k^2$. If a
	monomial $m$ does not contain submonomials $a_lb_k$, $b_la_k$, $a_la_k$,
	$b_lb_k$ for $l>k$, then the indices of variables in $m$ weakly
	increase from the left to right. We denote by $m_k$ a submonomial of
	$m$ that consists from the all variables with index $k$. Then $m=m_0 m_1
	\dots m_{l-1}$. Now if $m_k$ does not contain submonomals $a_k^2$,
	$b_k^2$, $b_ka_kb_ka_k$ then it is equal to one of the monomials
	$$
	e,\ a_k,\ b_k,\ a_kb_k,\ b_ka_k,\ a_kb_ka_k,\ b_ka_kb_k,\
	a_kb_ka_kb_k.
	$$
	Therefore there is no more then $8^l$ non-reducible monomials in $X^*_l$
	with respect to $S_l$. 
\end{proof}
\begin{corollary}
	\label{gb}
	Let $X=\bigcup_{l\le 0} X_l$. 
	The set $S=\bigcup_{l\ge 0}S_l$ is a reduced Gr\"obner basis of
	$\ker(\pi)$, where $\pi$ is the natural projection $\K\left\langle X^*
	\right\rangle\to \U_\K(sl^+_3)$. 
\end{corollary}
\begin{proof}
	It is clear that $S\subset \ker(\pi)$. Denote by $R$ the rewriting
	system $\left\{\, r(p) \,\middle|\, p\in S \right\}$. It is enough to show that   
	any critical pair $(w,r_1,r_2)$, with $r_1$, $r_2\in R$ is reducible.
	For a given critical pair $(w,r_1,r_2)$ there is an $l\ge 0$, such that
	all monomials in $w$, $r_1$, $r_2$ lie in $X_l$. By
	Proposition~\ref{gbl} the set $S_l$ is a Gr\"obner basis, therefore any
	critical pair $(w,r_1,r_2)$ with $w\in X_l^*$, $r_1$, $r_2\in
	R_l=\left\{\, r(p) \,\middle|\, p\in S_l \right\}$ is reducible.
\end{proof}

\section{First steps of minimal resolution for $n=3$ and $p=2$}
\label{firststeps}
We will consider the algebra $\U_\K(sl_3^+)$ as a graded algebra with the
grading induced by degree function $\deg(a_k)=\deg(b_k)=2^k$. Since the zero
component of $\U_\K(sl_3^+)$ is a field $\K$ and $\U_\K(sl_3^+)$ is a positively
graded all graded projective modules are shifts of the regular module
$\U_\K(sl_3^+)$. Suppose $P=\U_\K(sl^3_+)v$ with $\deg(v)=m$. We define
$Rad(R)= \bigoplus_{d\ge 1}\U_\K(sl_3^+)v$. This definition of radical can be
extended to the arbitrary projective module by additivity. It is well known that
a resolution of an $\U_\K(sl^+_3)$ module $M$ 
$$
\dots \to P_m \to \dots P_2 \stackrel{d_1}{\longrightarrow} P_1
\stackrel{d_0}{\longrightarrow} P_0 \to M \to 0
$$
is minimal up to step $m$ if and only if $\im(d_k) \subset Rad(P_k)$ for all
$0\le k\le m$. 

Now we will examine first four steps of the anick resolution of the trivial
module $\K$ over $\U_\K(sl^+_3)$. Then we modify it and get first three steps of
a minimal projective resolution. 

In our situation $T_0=\left\{\, a_k, b_k \,\middle|\, k\in \N_0 \right\}$,
and
$$T_1 =\left\{\, a_la_k, b_lb_k, a_lb_k, b_la_k, b_ka_kb_ka_kb_k, b_kb_k, a_ka_k
\,\middle|\, 0\le k< l \right\}.$$
We start with the computation of $d_1\colon P_1\to P_0$ in the anick resolution. 
Note that since we are working in characteristic two we can disregard sings. 
Suppose $l>k$, then
\begin{align*}
	d_1(.a_la_k) &= a_l.a_k + i_0d_0(a_l.a_k) = a_l.a_k + i_0(a_ka_l) \\
	&= a_l.a_k + a_k.a_l.
\end{align*}
Analogously $d_1(.b_lb_k)= b_l.b_k + b_k.b_l$.
Now
\begin{align*}
	d_1(.a_k^2) &= a_k.a_k + i_0d_0 (a_k.a_k) \\
	&= a_k.a_k + i_0(0) = a_k.a_k
\end{align*}
and analogously $d_1(.b_k^2)= b_k.b_k$.
\begin{align*}
	d_1(.a_lb_k) &= a_l.b_k + i_0d_0(a_l.b_k) \\
	&= a_l.b_k + i_0(a_kb_ka_k\dots a_{l-1}.e + b_ka_l.e) \\
	&= a_l.b_k + a_kb_ka_k\dots a_{l-2}. a_{l-1} + b_k.a_l.
\end{align*}
Analogously
\begin{align*}
	d_1(.b_la_k) &=  b_l.a_k + b_ka_kb_k\dots b_{l-2}. b_{l-1} + a_k.b_l.
\end{align*}
Now 
\begin{align*}
	d_1(.b_ka_kb_ka_k) &= b_ka_kb_k.a_k + i_0d_0(b_ka_kb_k.a_k)\\
	&= b_ka_kb_k.a_k + i_0(a_kb_ka_kb_k.e)\\
	&= b_ka_kb_k.a_k + a_kb_ka_k.b_k.
\end{align*}
It is readily seen that the image of $d_1$ is a subset of $Rad(P_0)$. 
Now we examine properties of $d_2\colon P_2\to P_1$ in the anick resolution. 
 In our case
\begin{align*}
	T_2 = & \left\{\, a_ma_la_k,\ b_mb_lb_k \,\middle|\,  m\ge l\ge k
	\right\} \\
	&\cup \left\{\, a_ma_lb_k,\ b_mb_la_k \,\middle|\,  m\ge l>k \right\}\\
	&\cup \left\{\, a_mb_la_k,\ b_ma_lb_k \,\middle|\,  m>l>k \right\}\\
	&\cup \left\{\, a_mb_lb_k,\ b_ma_la_k \,\middle|\,  m> l\ge k \right\}\\
	&\cup \left\{\, a_lb_ka_kb_ka_k,\ b_la_lb_la_lb_k \,\middle|\,  l>k
	\right\}\\
	&\cup \left\{\, b_lb_ka_kb_ka_k,\ b_la_lb_la_la_k \,\middle|\, l\ge k
	\right\}\\
	&\cup \left\{\, b_ka_kb_ka_kb_ka_k \,\middle|\, k\in \N_0 \right\}.
\end{align*}
We will not compute value of $d_2$ for all elements of $T_2$. Instead we will
show that for some elements of $T_2$ the image of $d_2$ lies in $Rad(P_1)$ and
for the rest of elements we compute $d_2$. Every element $v$  of $P_1$ can be
uniquely written as a sum
$$
\sum_{m\in M,\ t\in T_1} \lambda_{m,t}m.t,
$$
where $\lambda_{m,t}$ are elements of $\K$. Now $v\in Rad(P_1)$ if and only if
$\lambda_{e,t}=0$ for all $t\in T_1$. Therefore if $v$ is homogeneous and
$\deg(v)$ is not an element of the set $\left\{\, \deg(t) \,\middle|\, t\in T_1
\right\}$ then $v$ is an element of $Rad(P_1)$. Now
$$
D = \left\{\, \deg(t) \,\middle|\, t\in T_1 \right\} = \left\{\, 2^l+2^k \,\middle|\,
l\ge k\ge 0 \right\}.
$$
Note that $d_2$ preserves degree as our Gr\"obner basis is homogeneous.
Therefore if $t\in T_2$ and $\deg(t)\not\in D$, then $d_2(.t)\in Rad(P_1)$.
Now if $m>l>k$ then the degree of elements similar to $a_ma_la_k$ is
$2^m+2^l+2^k$ and it is not an element of $D$. Thus we have to consider only the
cases when at least two numbers $m$, $l$, $k$ are equal. 

Note  that $d_2(.b_ka_kb_ka_kb_ka_k)$, $d_2(.a_k^3)$, $d_2(.b_k^3)$,
$d_2(b_k^2a_kb_ka_k)$, $d_2(b_ka_kb_ka_k^2)$
are  linear combinations of monomials that
involve only variables with index $k$, since the linear span of such monomials is a
subalgebra of $\U_\K(sl^+_3)$. 
The set
$$
\left\{\, \deg(t) \,\middle|\, t\in T_1,\ t\in \left\{ a_k,b_k \right\}^* \right\}
$$
contains two numbers: $2^{k+1}$ and $2^{k+2}$. Now $\deg( (b_ka_k)^3)= 3\times
2^{k+1}$, $\deg(a_k^3)=\deg(b^3)= 3\times 2^k$,
$\deg(b_k^2a_kb_ka_k)=\deg(b_ka_kb_ka_k^2)=5\times 2^k$ are different from both
of them. Therefore $d_2(.(b_ka_k)^3)$, $d_2(.a_k^3)$, $d_2(.b_k^3)$,
$d_2(.b_k^2a_kb_ka_k)$, $d_2(.b_ka_kb_ka_k^2)$ lie in
the radical of $P_1$. 

Suppose $m>k$. Then the elements $a_m^2a_k$, $a_m^2b_k$, $b_m^2a_k$,
$b_m^2b_k$, $b_ma_mb_ma_ma_k$, and $b_ma_mb_ma_mb_k$ are elements of
the subalgebra generated by $X_m$. They all are of degree
$2^{m+1}+ 2^k$, but the set 
$$
\left\{\, \deg(t) \,\middle|\, t\in T_1,\ t\in X_m^* \right\} = 
 \left\{\, 2^l+2^k \,\middle|\,
m\ge l\ge k\ge 0 \right\}.
$$
does not contain $2^{m+1}+2^k$, therefore the value of $d_2$ for all
of the above
mentioned elements lies in the radical of $P_1$. 

It is left to compute $d_2$ for $a_ma_k^2$, $a_mb_k^2$, $b_ma_k^2$,
$a_mb_ka_kb_ka_k$, and $b_mb_ka_kb_ka_k$. 
Since $a_m$ and $a_k$
commute we get
\begin{align*}
d_2(.a_la_k^2) &= a_l.a_k^2 + a_k.a_la_k\\
\end{align*}
	and the similar formulae are valid for $d_2(.b_mb_k^2)$. 

Suppose  $m>k+1$. Then we have
\begin{align*}
	d_2(.a_mb_k^2) &= a_m.b_k^2 + i_1(a_mb_k.b_k) \\
	&= a_m.b_k^2 + i_1(b_ka_m.b_k + a_kb_ka_k\dots a_{m-1}.b_k) \\
	&= a_m.b_k^2 + b_k.a_mb_k + i_1(a_kb_ka_k\dots a_{m-1}.b_k +
	b_ka_kb_ka_k\dots a_{m-2}.a_{m-1})\\
	&= a_m.b_k^2 + b_k.a_mb_k+ a_kb_ka_k\dots a_{m-2}.a_{m-1}b_k
\end{align*}
and the similar formula for $d_2(.b_ma_k^2)$. Now suppose $m=k+1$. Then we get
\begin{align*}
	d_2(.a_{k+1}b_k^2) &= a_{k+1}.b_k^2 + i_1(a_{k+1}b_k.b_k)\\
	&= a_{k+1}.b_k^2+ i_1(b_ka_{k+1}.b_k + a_kb_ka_k.b_k)\\
	&= a_{k+1}.b_k^2 + b_k.a_{k+1}b_k + i_1(a_kb_ka_k.b_k + b_ka_kb_k.a_k)\\
	&= a_{k+1}.b_k^2 + b_k.a_{k+1}b_k + .b_ka_kb_ka_k
\end{align*}
and
\begin{align*}
	d_2(.b_{k+1}a_k^2) &= b_{k+1}.a_k^2 + a_k.b_{k+1}a_k + .b_ka_kb_ka_k.
\end{align*}
Now
\begin{align*}
	d_2(.a_mb_ka_kb_ka_k) &= a_m.b_ka_kb_ka_k + i_1(a_mb_ka_kb_k.a_k +
	a_ma_kb_ka_k.b_k)\\
	&= a_m.b_ka_kb_ka_k + i_1(b_ka_kb_ka_m.a_k+ a_kb_ka_ka_m.b_k)\\
	&= a_m.b_ka_kb_ka_k + b_ka_kb_k.a_ma_k +i_1(a_kb_ka_ka_m.b_k +
	b_ka_kb_ka_k.a_m)\\
&= a_m.b_ka_kb_ka_k + b_ka_kb_k.a_ma_k +i_1(a_kb_ka_ka_m.b_k +
	a_kb_ka_kb_k.a_m)\\
&= a_m.b_ka_kb_ka_k + b_ka_kb_k.a_ma_k +a_kb_ka_k.a_mb_k.
\end{align*}
\begin{align*}
 d_2(.b_mb_ka_kb_ka_k) & =  b_m.b_ka_kb_ka_k + i_1(b_mb_ka_kb_k.a_k +
 a_ma_kb_ka_k.b_k)\\
 & = b_m.b_ka_kb_ka_k + i_1(b_ka_kb_kb_m.a_k + a_kb_ka_ka_m.b_k)\\
 &= b_m.b_ka_kb_ka_k +b_ka_kb_k.b_ma_k + i_1(b_ka_kb_ka_k.b_m+
 a_kb_ka_ka_m.b_k)\\
&= b_m.b_ka_kb_ka_k +b_ka_kb_k.b_ma_k +i_1(a_kb_ka_kb_k.b_m+
 a_kb_ka_ka_m.b_k)\\
&= b_m.b_ka_kb_ka_k +b_ka_kb_k.b_ma_k +a_kb_ka_k.a_mb_k\ .
\end{align*}
Therefore $d_2(.a_mb_ka_kb_ka_k)$ and $d_2(.b_mb_ka_kb_ka_k)$ lie in
the radical of $P_1$. 

Let $T'_1=T_1\setminus \left\{ \,b_ka_kb_ka_k\,\middle|\, k\in \N_0
\right\}$ and $T'_2=T_2\setminus \left\{ \, b_{k+1}a_k^2\,\middle|\,
k\in \N_0 \right\}$. Define $P'_1$ and $P'_2$ to be the submodules of
$P_1$ and $P_2$ with $A$-bases $T'_1$ and $T'_2$ respectively. We define
$d'_1$ to be the restriction of $d_1$ on $P_1$. The differential
$d'_2\colon P'_2\to P'_1$ is defined as follows. Let $t\in T'_2$. Then
$$
d_2(.t) = f + \sum_{k=0}^{\infty} f_k.b_ka_kb_ka_k
$$
where $f\in P'_2$ and only finitely many $f_k$ are different from
$0$. Define
$$
d'_2(.t) = f + \sum_{k=0}f_k\left(b_{k+1}.a_k^2 + a_k.b_{k+1}a_k
\right).
$$
By the usual consideration we can see that the complex
$$
P'_2 \to P'_1 \to P'_0 \to P'_{-1} \to \K \to 0
$$
is exact. Moreover it is minimal up to the term $P'_1$. We get 
\begin{theorem}
 \label{main} Let us denote $\U_\K(sl^+_3)$ by $A$ and the free
 module over $A$, which is generated by an element of degree $s$, by $A[s]$. The
 trivial module $\K$ over $A$ has a minimal projective resolution of
 the form 
 $$
 \left( \bigoplus_{0\le k} A[2^{k+1}]^{\oplus 2} \right) \oplus
 \left( \bigoplus_{0\le k < l} A[2^l+2^k]^{\oplus 4} \right)\to
 \bigoplus_{k=0}^{\infty} A[2^k]^{\oplus 2} \to A \to \K\to 0.
 $$
\end{theorem}

\section{Conjectures}
\label{conjectures}

In this section we formulate several conjectures that were guessed from the excessive computer computations. 
We will consider the set of generators $X=\left\{\, e_{i,i+1}^{(p^k)} \,\middle|\, 1\le i\le n-1,\ k\in \N_0 \right\}$
of $\U_\K(sl^+_n)$, where $p$ is the characteristic of the field $\K$. To make formulae more readable we shall write $a_{ik}$ instead of $e_{i,i+1}^{(p^k)}$. 
We will assume the ordering $\le$ on $X$ defined by
$$
a_{11}<a_{21}<\dots<a_{n-1,1}<a_{12}<\dots <a_{n-1,2}<\dots  
$$
and on $\U_\K(sl^+_n)$ we consider deglex ordering, where $\deg(a_{ik})=p^k$.
\begin{conjecture}
\label{f}
	The map $X\to \U_\K(sl^+_n)$, $a_{ik}\mapsto a_{i,k+j}$ can be extended to a homomorphism of graded algebras $\U_\K(sl^+_n)\to \U_\K(sl^+_n)$. 
\end{conjecture}
Note that in the case of $n=3$ and $p=2$ the claim of the conjecture easily follows from Proposition~\ref{gb}. 

\begin{conjecture}
	Suppose $n=3$ and $p>2$. Denote $a_{1k}$ by $a_k$ and $a_{2k}$ by $b_k$.  Then the set 
	\begin{align*}
		&
		\left\{\, a_k^p,\ b_k^p,\ b_k^2a_k-2b_ka_kb_k+a_kb_k^2,\ b_k2a_k^2-2a_kb_ka_k+a_k^2b_k,\ (b_ka_k)^p-(a_kb_k)^p \,\middle|\, k\in\N_0 \right\}\\
		&
		\cup
		\left\{\, a_lb_k-b_ka_l-a_kb_ka_k^{p-1}\dots a_{l-1}^{p-1},\ b_la_k-a_kb_l-b_ka_kb_k^{p-1}\dots b_{l-1}^{p-1} \,\middle|\, l>k \right\}.\end{align*}
		is the Gr\"obner basis of $\U_\K(sl^+_4)$.
\end{conjecture}

Now we formulate a conjecture about Gr\"obner basis for the case $p=2$ and $n\ge 4$. For every sequence of integers $i=(i_1,\dots, i_l)$ we denote by $a_{i,k}$ the product
$$
a_{i_1,k}\dots a_{i_l,k}.
$$
We write $i..j$ for a sequence $(i,i+1,\dots,j)$ if $i<j$ and for a sequence $i,i-1,\dots,j$ if $i>j$.
We denote by $L_{l}$ the set of  all permutations $(i_1,\dots,i_l)$ of $(1,\dots,l)$ such that for every $1\le \sigma\le l-1$ either $i_{\sigma+1}<i_\sigma$ or $i_{\sigma+1}=i_\sigma+1$. Define
for every $m$ such that $m+l\le n$
$$
L_{l}[m] = \left\{\, (i_1+m,\dots,i_l+m) \,\middle|\, (i_1,\dots,i_l)\in L_{l} \right\}.
$$

\begin{conjecture}
	The set
	\begin{align*}
		&\left\{\, a_{ik}^2 \,\middle|\, 1\le i\le n-1,\ k\in \N_0 \right\}
		\cup \left\{\, \sum_{i\in L_l[m]} a_{ik}^2 \,\middle|\, k\in \N_0, l+m\le n \right\}\\
		&\cup\left\{\, [a_{i+m-1,k},a_{i..i+m,k}]+[a_{i+m-1,k},a_{i+m..i,k}] \,\middle|\, 1\le i\le n-m-1,\ k\in \N_0 \right\}\\
		&\cup \left\{\, a_lb_k+b_ka_l + a_kb_ka_k\dots a_{l-1},\ b_la_k+a_kb_l+b_ka_kb_k\dots b_{l-1} \,\middle|\,  \right\}.
	\end{align*}
is a Gr\"obner basis of $\U_\K(sl^+_n)$.
\end{conjecture}
Note that this conjecture has a very simple proof module Conjecture~\ref{f}.
\bibliography{groebnersl2}

\providecommand{\bysame}{\leavevmode\hbox to3em{\hrulefill}\thinspace}
\providecommand{\MR}{\relax\ifhmode\unskip\space\fi MR }
\providecommand{\MRhref}[2]{%
  \href{http://www.ams.org/mathscinet-getitem?mr=#1}{#2}
}
\providecommand{\href}[2]{#2}
\begin{thebibliography}{1}

\bibitem{anick}
D.J. Anick, \emph{On the homology of associative algebras}, Trans. Amer. Math.
  Soc. \textbf{296} (1986), no.~2, 641--659.

\bibitem{kostant}
Bertram Kostant, \emph{Groups over {$Z$}}, Algebraic {G}roups and
  {D}iscontinuous {S}ubgroups ({P}roc. {S}ympos. {P}ure {M}ath., {B}oulder,
  {C}olo., 1965), Amer. Math. Soc., Providence, R.I., 1966, pp.~90--98.
  \MR{MR0207713 (34 \#7528)}

\bibitem{apsme}
Ana~Paula Santana and Ivan Yudin, \emph{The {K}ostant form of
  {$\mathfrak{U}(sl_n^+)$} and the {B}orel sublagebra of the {S}chur algebra
  {$S(n,r)$}}, \path|arXiv:0803.4382|.

\end{thebibliography}
\bibliographystyle{amsplain}

\end{document}